\documentclass[11pt,reqno]{amsart}
\usepackage{a4wide}

\usepackage{amssymb, amsbsy, amsfonts}
\usepackage{amsmath, amsthm, mathrsfs}
\usepackage{latexsym,float,color}
\usepackage{etex}
\usepackage{eucal}
\usepackage{url}

\numberwithin{equation}{section}

\usepackage{tikz}
\usetikzlibrary{matrix}

\definecolor{blaugrau}{rgb}{0.796887, 0.789075, 0.871107}

\newcounter{linectr}

\allowdisplaybreaks[4]

\newcounter{mmacnt}
\def\restartmma{\setcounter{mmacnt}{0}}
\restartmma \catcode`|=\active
\def|#1|{\mathrm{#1}}
\catcode`|=12
\newenvironment{mma}{
 \par
 \catcode`|=\active
 \parskip=2pt\parindent=0pt 
 \small
 \def\In##1\\{%
   \def\linebreak{\hfill\break\null\qquad}%
   \refstepcounter{mmacnt}
   \hangindent=2.5em\hangafter=0
   \leavevmode
   \llap{\tiny\sffamily In[\arabic{mmacnt}]:=\kern.5em}%
   \mathversion{bold}\scriptsize$\tt\bf\displaystyle##1$\normalsize
   \mathversion{normal}\par
 }%
 \def\Print##1\\{%
   \def\linebreak{\hfill\break}%
   \hangindent=2.5em\hangafter=0
   \leavevmode\scriptsize ##1\par}%
 \def\Out##1\\{%
   \vspace*{-0.2cm}\def\linebreak{$\hfill\break\null\hfill$}%
   \kern\abovedisplayskip\par
   \hangindent=2.5em\hangafter=0
   \leavevmode
   \llap{\tiny\sffamily Out[\arabic{mmacnt}]=\kern.5em}
   \scriptsize$\displaystyle\tt##1$\normalsize\hfill\null\par
   \kern\belowdisplayskip\vspace*{-0.3cm}
 }%
 \def\Warning##1##2\\{%
   \def\linebreak{\hfill\break}%
   \hangindent=2.5em\hangafter=0
   \leavevmode
   {\scriptsize##1 : ##2}\par}%
}{%
 \par\smallskip
}

\newcommand{\LoadP}[1]{\fcolorbox{black}{blaugrau}{
\begin{minipage}[t]{14cm}
\footnotesize #1
\end{minipage}}}

\def\<#1>{\langle#1\rangle}
\let\set\mathbb

\def\NN{{\set Z}_{\geq0}}

\newdimen\listablecorrection
\theoremstyle{definition}
\newtheorem{theorem}{Theorem}[section]

\newtheorem{definition}[theorem]{Definition}

\title[A Motivating Challenge Solved by Computer Algebra]{The Absent-Minded Passengers Problem: A Motivating Challenge Solved by Computer Algebra}

\author{Carsten Schneider}
\address[Carsten Schneider]{Johannes Kepler University Linz\\
Research Institute for Symbolic Computation\\
A-4040 Linz, Altenberger Stra{\ss}e 69, Austria}
\email{Carsten.Schneider@risc.jku.at}

\thanks{Supported by the Austrian Science Fund (FWF) grant SFB F50 (F5009-N15).}

\begin{document}

\begin{abstract}
In~\cite{EZ:20} an exciting case study has been initiated in which experimental mathematics and symbolic computation are utilized to discover new properties concerning the so-called Absent-Minded Passengers Problem. Based on these results, Doron Zeilberger raised some challenging tasks to gain further probabilistic insight. In this note we report on this enterprise. In particular, we demonstrate how the computer algebra packages of RISC can be used to carry out the underlying heavy calculations. 
\end{abstract}

\maketitle

\section{A challenging Email}

On January 22, 2020 I received the following email by Doron Zeilberger:

\medskip

\begin{verbatim}
Dear Carsten,

I (and Shalosh) just posted a paper

https://arxiv.org/abs/2001.06839

with a challenge to you (see the middle of page 4)

Can you (and Sigma) extend theorem 5 of that paper
to the general case with k absent-minded passengers?

These expressions should be expressions in n,k,
the harmonic numbers H_{n-1}, H_{k-1} and their
generalizations (the partial sum , to n-1, k-1, respectively of
zeta(n) and zeta(k)).

If you and Sigma can do the fourth moment, and derive
the asymptotic in n (with a fixed but arbitrary k), I will
donate $100$ to the OEIS in your honor.

...

Best wishes,
Doron


\end{verbatim}

\medskip

When I received this email, I was thrilled:
First, if one gets such an email -- in particular from Doron Zeilberger, one is automatically eager to solve it.
And second, I highly appreciate ``The On-Line Encyclopedia of Integer Sequences'' (OEIS) at \texttt{http://www.oeis.org} and supporting it by a donation of Doron Zeilberger gave an extra strong motivation.

Summarizing, the above email provoked various heavy calculations by means of computer algebra that will be described in the following.

\section{The underlying problem and symbolic summation}

The combinatorial problem of absent-minded passengers can be introduced as follows.

\medskip

\begin{definition}\label{Def:AbsentMindedPassenger}
Consider a plane with $n\geq2$ seats and suppose that $n$ passengers enter the plane step-wise taking their seats. In addition, suppose that the first $k\geq1$ passengers are absent-minded, i.e., they lost their seat ticket and take a seat uniformly at random. The remaining $n-k$ not absent-minded (but shy) passengers take their dedicated seats (as given in the plane ticket) if it is still free; otherwise, they choose uniformly at random one of the still available free seats.  
For $0\leq r\leq n$, let $p_{n,k,r}$ be the probability that exactly $r$ passengers sit in the wrong seat 
and let $X_n$ be the random variable for ``the number of passengers sitting in the wrong seat''.
Then the expected value for the passengers sitting in the wrong seat is
\begin{align*}
E(X_n)&=\sum_{r=1}^{n}r\,p_{n,k,r}\\[-0.2cm]
\intertext{and its variance is} 
V(X_n)&=E(X_n^2)-E(X_n)^2=\sum_{r=1}^{n}r^2\,p_{n,k,r}-E(X_n)^2.
\end{align*}
\end{definition}

The situation of one absent-minded passenger ($k=1$) has been considered in~\cite{Winkler:04,Bolobas:06}
and has been explored further in~\cite{HL:19} for the general case $k\geq1$.
Among other fascinating results, closed forms for $E(X_n)$ and $V(X_n)$ have been obtained in~\cite{HL:19} by skillful combinatorial arguments. More precisely, the definite sum representations 
\begin{align}
\label{Equ:ENotSimplified}
E(X_n)&=\frac{k (n-1)}{n}
+
\sum_{i=1}^{-k
	+n
} \frac{k}{1
	-i
	+n
}\\[-0.2cm]
\intertext{and}
\label{Equ:VNotSimplified}
V(X_n)&=\frac{k (n-1)}{n^2}+\sum_{i=1}^{-k
	+n
} \frac{(1
	-i
	-k
	+n
	) \big(
	1
	-\frac{1
		-i
		-k
		+n
	}{1
		-i
		+n
	}
	\big)}{1
	-i
	+n
}\\[-0.2cm]
&\hspace*{0.5cm}+2 \left(
\frac{(k-1) k}{2 (n-1) n^2}
+
\sum_{i=1}^k 
\sum_{j=1}^{-k
	+n
} \frac{\frac{1
		-j
		-k
		+n
	}{-j
		+n
	}
	-\frac{1
		-j
		-k
		+n
	}{1
		-j
		+n
	}
}{n}
\right)\nonumber
\end{align}
have been derived and simplified to 
\begin{align}
\label{Equ:ESimplified}
E(X_n)&=\frac{k (n-1)}{n}-k S_1(k)+k S_1(n)\\
\intertext{and}
\label{Equ:VSimplified}
V(X_n)&=\frac{
	2 k
	-k^2
	-2 n
	-2 k n
	+2 k^2 n
	+2 n^2
	-k n^2
}{(n-1) n^2}\\
\nonumber
&-\frac{k (2+n) S_1(k)}{n}
+\frac{k (2+n) S_1(n)}{n}
+k^2 S_2(k)
-k^2 S_2(n)
\end{align}
in terms of the harmonic numbers  
$$S_o(n)=\sum_{i=1}^n\frac1{i^o}$$
of order $o\geq1$ (often they are also denoted by $H^{(o)}_n=S_o(n)$ with the special case $H_n=S_1(n)$). In the following we will prefer to write such expressions in terms of the modified harmonic numbers 
$$\bar{S}_o(n)=S_o(n-1)=\sum_{i=1}^{n-1}\frac1{i^o}$$
yielding, e.g., the more compact expressions
\begin{align}
\label{Equ:EValue}
E(X_n)&=-1
+k
-k \bar{S}_1(k)
+k \bar{S}_1(n)\\
\intertext{and}
\label{Equ:VValue}
V(X_n)&=\frac{(k-1) k}{(n-1) n}
-\frac{k (2+n) \bar{S}_1(k)}{n}
+\frac{k (2+n) \bar{S}_1(n)}{n}
+k^2 \bar{S}_2(k)
-k^2 \bar{S}_2(n).
\end{align}
While the simplification from~\eqref{Equ:ENotSimplified} to~\eqref{Equ:ESimplified} is straightforward, more work has to be carried out to derive~\eqref{Equ:VSimplified} from~\eqref{Equ:VNotSimplified}. In~\cite{HL:19} further details are suppressed how these simplifications have been obtained. Here I want to point out that such (usually painful) classical manipulations are meanwhile obsolete. For instance, by loading in the computer algebra packages 
\begin{mma}
	\In << Sigma.m \\
	\Print \LoadP{Sigma - A summation package by Carsten Schneider
		\copyright\ RISC-JKU}\\
	\In << HarmonicSums.m\\
	\Print \LoadP{HarmonicSums by Jakob Ablinger
	\copyright\ RISC-JKU}\\
	\In << EvaluateMultiSums.m\\
\Print \LoadP{EvaluateMultiSums by Carsten Scneider
	\copyright\ RISC-JKU}\\
\end{mma}

\medskip

\noindent into Mathematica the multi-sum expression for the variance

\medskip

\begin{mma}
\In V=\frac{k (n-1)}{n^2}+
\sum_{i=1}^{-k
	+n
} \frac{(1
	-i
	-k
	+n
	) \big(
	1
	-\frac{1
		-i
		-k
		+n
	}{1
		-i
		+n
	}
	\big)}{1
	-i
	+n
}\newline
\hspace*{0.5cm}+2 \left(
\frac{(k-1) k}{2 (n-1) n^2}
+
\sum_{i=1}^k 
\sum_{j=1}^{-k
	+n
} \frac{\frac{1
		-j
		-k
		+n
	}{-j
		+n
	}
	-\frac{1
		-j
		-k
		+n
	}{1
		-j
		+n
	}
}{n}
\right)
;\\
\end{mma}

\medskip

\noindent can be simplified to a closed form within seconds by executing the command\footnote{Within the \texttt{Sigma} package (and later the \texttt{HarmonicSums} package) $S_o(n)$ is denoted by $\texttt{S[o,n]}$.}

\medskip

\begin{mma}
	\In EvaluateMultiSum[V,\{\},\{k,n\},\{1,2\},\{n,Infinity\}]\\
	\Out \frac{
	2 k
	-k^2
	-2 n
	-2 k n
	+2 k^2 n
	+2 n^2
	-k n^2
	}{(n-1) n^2}
	-\frac{k (2+n) S[1,k]}{n}
	+\frac{k (2+n) S[1,n]}{n}
	+k^2 S[2,k]
	-k^2 S[2,n]\\
\end{mma}

\medskip

\noindent\textit{Remark.} Internally, the command \texttt{EvaluateMultiSums} of the package \texttt{EvaluateMultiSums.m} \cite{Schneider:13a} uses the summation paradigms of telescoping, creative telescoping and recurrence solving of the summation package \texttt{Sigma.m}~\cite{Schneider:07a}; the underlying algorithms generalize the hypergeometric case~\cite{AequalB} to the class of indefinite nested sums and products in the setting of difference fields and rings~\cite{Karr:81,DR1,DR3}. In addition, the calculations are supported by special function algorithms of the package \texttt{HarmonicSums.m}~\cite{ABS:11,ABS:13}. We note that the user is completely freed from applying the summation tools explicitly. However, all the calculation steps are equipped with proof certificates (based on the creative telescoping paradigm introduced in~\cite{Zeilberger:91}). Thus if necessary, a rigorous correctness proof can be extracted.

\section{The generating function approach for large moments}

As elaborated in detail in~\cite{EZ:20} the expectation, the variance and higher moments can be elegantly described with the generating function approach. Namely, consider the generating function
$$f_n^{(k)}(w)=\sum_{r=0}^{n}p_{n,k,r}\,w^r$$
of the probability $p_{n,k,r}$ introduced in Definition~\ref{Def:AbsentMindedPassenger}. Then 
the expectation and variance (see also Definition~\ref{Def:AbsentMindedPassenger}) 
can be straightforwardly connected to the generating function via
$$E(X_n)=\frac{d}{dw}f_n^{(k)}(w)|_{w=1}$$
and
$$V(X_n)=\big(w\tfrac{d}{dw}\big)^2f_n^{(k)}(w)|_{w=1}-E(X_n)^2.$$
More generally, define
$$M_0(n,k)=1$$
and
\begin{equation}\label{Equ:ComputeMi}
M_l(n,k)=\big(w\tfrac{d}{dw}\big)^lf_{n}^{(k)}(w)|_{w=1}
\end{equation}
for $l\geq1$; note that $M_1(n,k)=E(X_n)$. Then one can define the $l$-th moment by
\begin{equation}\label{Equ:Formulaml}
m_l(n,k)=\sum_{i=0}^l\binom{l}{i}(-1)^{l-i}M_i(n,k)M_{1}^{l-i}(n,k).
\end{equation}

In the above email but also in~\cite{EZ:20} the task was assigned to compute (besides the known moments $m_1(n,k)=0$ and $m_2(n,k)=V(X_n)$) further moments $m_l(n,k)$ (at least for $l=3,4$).

\medskip 

The combinatorial approach of~\cite{HL:19} to derive multi-sum expressions of $m_l(n,k)$ for larger $l$ seems hopeless. Even though it would have been pure fun to challenge \texttt{Sigma.m}, similarly as in in~\cite{PSW:11}, with more complicated sums than~\eqref{Equ:ESimplified}.
Another exciting approach has been carried out in~\cite[Theorem~5]{EZ:20} by guessing the moments $m_l(n,k)$ for the special case $k=1$ and $l=2,\dots,8$.  
As it turns out (and proposed by Doron Zeilberger in his email), this attempt can be pushed forward by the following powerful formula given in~\cite[Theorem~3]{EZ:20}:
\begin{equation}\label{Equ:Doronfn}
\begin{split}
f_n^{(k)}(w)&=\frac1{n}\sum_{r=0}^kr!\binom{k}{r}w^r(1-w)^{k-r}\prod_{i=0}^{n-k-1}(r\,w+i+1)\\
&=\frac{1}{n!}\sum_{r=0}^k(1+r\,w)r!\binom{k}{r}w^r(1-w)^{k-r}(2+r\,w)_{n-k-1};
\end{split}
\end{equation}
here $(x)_k=x(x+1)\dots(x+k-1)$ denotes the Pochhammer symbol. As a consequence, one can calculate straightforwardly any value of $M_l(n,k)$ and thus of $m_l(n,k)$ for particularly chosen $n$ and $k$. E.g., we can compute
\medskip
\begin{mma}
\In F=\frac{1}{n!}\sum_{r=0}^k(1+r\,w)r!\binom{k}{r}w^r(1-w)^{k-r}(2+r\,w)_{n-k-1}/.\{n\to10,k\to2\};\\
\In MSpec0=F/.w\to1\\
\Out 1\\
\vspace*{0.2cm}
\In MSpec1=w*D[F,w]/.w\to1\\
\Out \frac{5869}{1260}\\
\vspace*{0.2cm}
\In MSpec2=w*D[w*D[F,w],w]/.w\to1\\
\Out \frac{50293}{2100}\\
\vspace*{0.2cm}
\In MSpec3=w*D[w*D[w*D[F,w],w],w]/.w\to1\\
\Out \frac{9966821}{75600}\\
\end{mma}

\medskip

\noindent and activating the formula~\eqref{Equ:Formulaml} for $n=10$, $k=2$ and $l=3$ yields the third moment

\medskip

\begin{mma}
\In mSpec3=-MSpec0*MSpec1^3+3*MSpec1^3-3*MSpec2*MSpec1+MSpec3\\
\Out -\frac{702653939}{1000188000}\\
\end{mma}
\medskip

\noindent But even more is possible with the hypergeometric sum representation~\eqref{Equ:Doronfn}. As already proposed in~\cite{EZ:20} one can derive easily the closed form expressions of $M_l(n,k)$ for symbolic $n$ and $k$. For this task we observe that
\begin{align}
M_l(n,k)=&\big(w\tfrac{d}{dw}\big)^lf_{n}^{(k)}(w)|_{w=1}\nonumber\\
=&\big(w\tfrac{d}{dw}\big)^l\left[\frac{1}{n!}\sum_{r=0}^k(1+r\,w)r!\binom{k}{r}w^r(1-w)^{k-r}(2+r\,w)_{n-k-1}\right]\Big|_{w=1}\nonumber\\
\label{Equ:DoronGenFu}
=&\big(w\tfrac{d}{dw}\big)^l\left[\frac{1}{n!}\sum_{r=\max(0,k-l)}^k(1+r\,w)r!\binom{k}{r}w^r(1-w)^{k-r}(2+r\,w)_{n-k-1}\right]\Big|_{w=1},
\end{align}
i.e., for any $l\in\NN$ at most $l+1$ summands contribute and the remaining summands vanish with the evaluation $w=1$.
Moreover the differentiation of the arising building blocks in~\eqref{Equ:DoronGenFu} can be carried out easily. For instance if we differentiate $(2+r\,w)_{n-k-1}$ twice w.r.t.\ $w$ by using the Mathematica-command \texttt{D} and the \texttt{Sigma}-commands \texttt{ToSigma} and \texttt{CollectProdSum} we get

\begin{mma}
\In D[D[(2+r\,w)_{n-k-1},w],w]//ToSigma//CollectProdSum\\
\Out \Bigg(
\frac{2 r^2}{(1
	+r w
	)^2}
+\frac{2 r^2 S[1,{r w}]}{1
	+r w
}
+r^2 S[1,{r w}]^2
-\frac{2 r^2 S[1,{-k+n+r w}]}{1
	+r w
}
-2 r^2 S[1,{r w}] S[1,{-k+n+r w}]
+r^2 S[1,{-k+n+r w}]^2\newline
\hspace*{0.5cm}+r^2 S[2,{r w}]
-r^2 S[2,{-k+n+r w}]
\Bigg) (2
+r w
)_{-1
	-k
	+n
}\\
\end{mma}

\medskip

In this way, one can compute $M_1(n,k)(=E(X_n))$ and rediscovers~\eqref{Equ:EValue}. Similarly, one gets, e.g.,
\begin{equation}\label{Equ:M2}
\begin{split}
M_2(n,k)&=(k-1)^2
+\frac{\big(
	-2 k
	+k n
	-2 k^2 n
	\big) \bar{S}_1(k)}{n}
+k^2 \big(
\bar{S}_1(k)\big)^2\\
&+\frac{\big(
	2 k
	-k n
	+2 k^2 n
	\big) \bar{S}_1(n)}{n}
-2 k^2 \bar{S}_1(k) \bar{S}_1(n)
+k^2 \big(
\bar{S}_1(n)\big)^2
+k^2 \bar{S}_2(k)
-k^2 \bar{S}_2(n)
\end{split}
\end{equation}
and
\begin{equation}\label{Equ:M3}
\begin{split}
M_3(n,k)&=
2 k^3 \bar{S}_3(n)-2 k^3 \bar{S}_3(k)-\frac{k \left(6 k-n^2-5 n\right) \bar{S}_1(n)}{(n-1)
	n}-\frac{3 k \bar{S}_1(k)^2}{n}-\frac{3 k \bar{S}_1(n)^2}{n}\\
&+\bar{S}_1(k) \left(\frac{6 k \bar{S}_1(n)}{n}+\frac{k \left(6 k-n^2-5 n\right)}{(n-1)
	n}\right)+\frac{3 k (k n+2 k-1)
	\bar{S}_2(k)}{n}\\
&-\frac{3 k (k n+2 k-1) \bar{S}_2(n)}{n}-\frac{(k-2) (k-1) k}{(n-2) (n-1) n}+\frac{3 (k-2) k (n-k)}{(n-1)^2 n}\delta(k-2)
\end{split}
\end{equation}
with
$$\delta(x)=\begin{cases} 1 &\text{ if } x\geq0\\ 0&\text{ if } x<0.\end{cases}$$
Note that in $M_3(n,k)$
the factor $\frac{3 (k-2) k (n-k)}{(n-1)^2 n}$ does not contribute for $k=1,2$ ($k=2$ would vanish anyway); this comes from the fact that in the summation of~\eqref{Equ:DoronGenFu} the lower bound $k-l$ with $l=3$ should be non-negative.  

In a nutshell, one can calculate sufficiently many $M_l(n,k)$ and using the formula~\eqref{Equ:Formulaml} one gets the desired moments $m_l(n,k)$. With $M_1(n,k)=E(X_n)$ given in~\eqref{Equ:EValue},~\eqref{Equ:M2} and~\eqref{Equ:M3} one can reproduce $m_2(n,k)$ as given in~\eqref{Equ:VValue} and can compute in addition
\small
\begin{align*}
m_3(n,k)=&-\frac{(k-2) (k-1) k}{(n-2) (n-1) n}+2 k^3 \bar{S}_3(n)
+\frac{\big(
	6 k^2
	-5 k n
	-k n^2
	\big) \bar{S}_1(k)}{(n-1) n}
-\frac{3 k \big(
	\bar{S}_1(k)\big)^2}{n}\\
&+\frac{\big(
	-6 k^2
	+5 k n
	+k n^2
	\big) \bar{S}_1(n)}{(n-1) n}
+\frac{6 k \bar{S}_1(k) \bar{S}_1(n)}{n}
-\frac{3 k \big(
	\bar{S}_1(n)\big)^2}{n}+\frac{3 \big(
	-k
	+2 k^2
	+k^2 n
	\big) \bar{S}_2(k)}{n}\\
&-2 k^3 \bar{S}_3(k)
-\frac{3 \big(
	-k
	+2 k^2
	+k^2 n
	\big) \bar{S}_2(n)}{n}
+\frac{3 (k-2) k (n-k)}{(n-1)^2 n}\delta(k-2).\\
\intertext{\normalsize Similarly, one obtains, e.g.,}
m_4(n,k)=&
\frac{(k-3) (k-2) (k-1) k}{(n-3) (n-2) (n-1) n}
+\frac{
	6 k
	-12 k^2
	+30 k^3
	-6 k^4
	+18 k n
	-29 k^2 n
	-7 k^2 n^2
	}{(n-1) n} \bar{S}_2(n)\\
&+\frac{\bar{S}_2(k)}{(n-1) n} \big(
-6 k
+12 k^2
-30 k^3
+6 k^4
-18 k n
+29 k^2 n
+7 k^2 n^2
\big)\\
&+\frac{\bar{S}_1(k)}{(n-2) (n-1) n} \big(
-52 k
+36 k^2
-12 k^3
+40 k n
-11 k n^2
-k n^3
\big)\\
&+\frac{\bar{S}_1(n)}{(n-2) (n-1) n} \big(
52 k
-36 k^2
+12 k^3
-40 k n
+11 k n^2
+k n^3
\big)\\
&+\frac{3 \big(
	-2 k
	+4 k^2
	-6 k n
	+3 k^2 n
	+k^2 n^2
	\big)
	\big(\bar{S}_1(k)\big)^2}{(n-1) n}
-\frac{4 k \big(
	\bar{S}_1(k)\big)^3}{n}+3 k^4 \big(
\bar{S}_2(n)\big)^2\\
&-\frac{6 \big(
	-2 k
	+4 k^2
	-6 k n
	+3 k^2 n
	+k^2 n^2
	\big) \bar{S}_1(k) \bar{S}_1(n)}{(n-1) n}
+\frac{12 k \big(
	\bar{S}_1(k)\big)^2 \bar{S}_1(n)}{n}\\
&+\frac{3 \big(
	-2 k
	+4 k^2
	-6 k n
	+3 k^2 n
	+k^2 n^2
	\big)
	\big(\bar{S}_1(n)\big)^2}{(n-1) n}
-\frac{12 k \bar{S}_1(k) \big(
	\bar{S}_1(n)\big)^2}{n}
+\frac{4 k \big(
	\bar{S}_1(n)\big)^3}{n}\\
&-\frac{6 \big(
	2 k
	-4 k^2
	+2 k^3
	+k^3 n
	\big) \bar{S}_1(k) \bar{S}_2(k)}{n}
+\frac{6 \big(
	2 k
	-4 k^2
	+2 k^3
	+k^3 n
	\big) \bar{S}_1(n) \bar{S}_2(k)}{n}\\
&-\frac{4 \big(
	2 k
	-6 k^2
	+6 k^3
	+3 k^3 n
	\big) \bar{S}_3(k)}{n}
+6 k^4 \bar{S}_4(k)
+\frac{6 \big(
	2 k
	-4 k^2
	+2 k^3
	+k^3 n
	\big) \bar{S}_1(k) \bar{S}_2(n)}{n}\\
&-\frac{6 \big(
	2 k
	-4 k^2
	+2 k^3
	+k^3 n
	\big) \bar{S}_1(n) \bar{S}_2(n)}{n}
-6 k^4 \bar{S}_2(k) \bar{S}_2(n)\\
&+\frac{4 \big(
	2 k
	-6 k^2
	+6 k^3
	+3 k^3 n
	\big) \bar{S}_3(n)}{n}
-6 k^4 \bar{S}_4(n)-\frac{2}{(n-2)^2 (n-1)^2 n} \big(
-4 k
+36 k^2\\
&-30 k^3
+6 k^4
-24 k n
+k^2 n
+15 k^3 n
-4 k^4 n+16 k n^2
-15 k^2 n^2
+3 k^3 n^2
\big)\delta(k-3)\\
&+\Big[
-\frac{6 \bar{S}_1(k)}{(n-1)^2 n} \big(
-2 k
-5 k^2
+3 k^3
+10 k n
-7 k^2 n
+k^3 n
\big)\\
&+\frac{6 \bar{S}_1(n)}{(n-1)^2 n} \big(
-2 k
-5 k^2
+3 k^3
+10 k n
-7 k^2 n
+k^3 n
\big)+\frac{1}{(n-1)^3 n} \big(
-k
+37 k^2\\
&-42 k^3
+12 k^4
-34 k n
+52 k^2 n
-18 k^3 n
-13 k n^2
+7 k^2 n^2
\big)
\Big]\delta(k-2).\end{align*}
\normalsize
For instance, specializing $m_4(n,k)$ to $k=1$ gives
\begin{align*}
m_4(n,1)&=\big(
\frac{14+n}{n}
-6 \bar{S}_2({n})
\big) \bar{S}_1({n})
+\frac{3 (n-2) \bar{S}_1({n})^2}{n}
+\frac{4 \bar{S}_1({n})^3}{n}\\
&
-\frac{(18+7 n) \bar{S}_2({n})}{n}+3 \bar{S}_2({n})^2
+\frac{4 (2+3 n) \bar{S}_3({n})}{n}
-6 \bar{S}_4({n})
\end{align*}
for one absent-minded passenger (as derived in~\cite[Theorem~5]{EZ:20}) and 
\begin{align*}
m_4(n,2)&=\big(
\frac{2 \big(
	-74+13 n+13 n^2\big)}{(n-1) n}
-\frac{24 (1+2 n) \bar{S}_2({n})}{n}
\big) \bar{S}_1({n})+48 \bar{S}_2({n})^2\\
&+\frac{12 \big(
	5-2 n+n^2\big) \bar{S}_1({n})^2}{(n-1) n}
+\frac{8 \bar{S}_1({n})^3}{n}
-\frac{4 (21+19 n) \bar{S}_2({n})}{n}\\
&+\frac{16 (7+6 n) \bar{S}_3({n})}{n}
-96 \bar{S}_4({n})+\frac{2 \big(
	51-45 n+19 n^2\big)}{(n-1) n}
\end{align*}
for two absent-minded passengers.

Following this strategy we succeeded in computing the moments $m_l(n,k)$ up to $l=15$ on a machine with 12 cores and 1.5TB memory, but failed to proceed due to time and space limitations. Luckily, another great trick from~\cite{Zeil:09} enabled us to compute even more moments. Namely, the calculation of
\begin{equation}\label{Equ:FactMoments}
\begin{split}
\bar{M}_l(n,k)&=\big(\tfrac{d}{dw}\big)^lf_{n}^{(k)}(w)|_{w=1}\\
=&\big(\tfrac{d}{dw}\big)^l\left[\frac{1}{n!}\sum_{r=\max(0,k-l)}^k(1+r\,w)r!\binom{k}{r}w^r(1-w)^{k-r}(2+r\,w)_{n-k-1}\right]\Big|_{w=1}
\end{split}
\end{equation}
turns out to be much faster (and less memory consuming) than the calculation of~\eqref{Equ:DoronGenFu}. Given these so-called exponential moments, one gets back the ordinary moments $M_l(n,k)$ by using the following formula from~\cite{Zeil:09}:
\begin{equation}\label{Equ:TransformExpToOrd}
M_l(n,k)=\sum_{r=1}^lS(l,r)\bar{M}_l(n,k)
\end{equation}
where $S(l,r)$ denotes the Stirling numbers of the second kind. In summary, we carried out the following steps:
\begin{enumerate}
\item[Step 1:] Calculate the exponential moments $\bar{M}_l(n,k)$ with the formula~\eqref{Equ:FactMoments}.
\item[Step 2:] Calculate the ordinary moments $M_l(n,k)$ with the formula~\eqref{Equ:TransformExpToOrd}.
\item[Step 3:] Finally, calculate the desired moments $m_l(n,k)$ with the formula~\eqref{Equ:Formulaml}.
\end{enumerate}
Using this approach, we succeeded in calculating generously\footnote{Originally only $l\leq 4$ was required in the challenge.} the moments $m_l(n,k)$ up to $l=21$ in a decent amount of time and memory. The moments up to order 16 are available online at 
\label{Page:Links}
\begin{center}
\url{https://www.risc.jku.at/people/cschneid/data/AMPassenger.tar.gz}.
\end{center}
Furthermore, a Mathematica notebook is provided at
\begin{center}
\url{https://www.risc.jku.at/people/cschneid/data/AMPassenger.nb}
\end{center}
that illustrates the main calculations (e.g., up to $l=7$). For higher moments the calculations have been distributed (parallelized) neatly using about 10 Mathematica subkernels. Since the code is rather technical and not easily digestible, it is suppressed within the attached Mathematica notebook.
The corresponding timings for the various steps, the total time and the size of the moments are summarized in the table of Figure~\ref{Fig:Timing}. 

\begin{figure}
	\begin{tabular}{|c||c|c|c|c|c|}
		\hline
		l&time for step 1 &time for step 2 &time for step 3&total time &size of $m_l(n,k)$\\
		\hline
		\hline
		1 & 0.16 s& 0.016 s & 0.96 s& 1.14 s  &
		88 B \\
		2 & 1.86 s& 0.025 s& 0.068 s& 1.95 s&
		1536 B\\
		3 & 3.06 s& 0.067 s& 0.182 s& 3.307 s&
		6 KB\\
		4 & 4.61 s& 0.15 s& 0.48 s& 5.24 s&
		25 KB\\
		5 & 6.66 s& 0.33 s& 1.18 s& 8.18 s&
		73 KB\\
		6 & 9.55 s& 0.75 s& 2.81 s& 13.12 s&
		210 KB\\
		7 & 13.81 s& 1.70 s& 6.51 s& 22.02 s&
		531 KB \\
		8 & 21.53 s & 3.73 s & 15.14 s& 40.40 s&
		1.3 MB \\
		9 & 35.79 s& 8.03 s& 34.50 s& 78.33 s&
		2.9 MB \\
		10 & 63 s& 18 s& 75 s& 156 s&
		6.4 MB \\
		11 & 115 s& 37 s& 160 s& 313 s&
		13 MB \\
		12 & 222 s& 77 s& 335 s& 634 s&
		27 MB \\
		13 & 454 s& 233 s& 678 s& 1364 s&
		52 MB \\
		14 & 1101 s& 579 s& 1344 s& 3024 s&
		98 MB \\
		15 & 2559 s& 1063 s& 2611 s& 6233 s&
		180 MB \\
		16 & 5249 s& 2380 s& 4995 s& 12625 s&
		326 MB \\
		17 & 11510 s& 4238 s& 9521 s& 25270 s&
		573 MB\\
		18 & 22357 s& 8807 s& 17669 s& 48834 s&
		993 MB\\
		19 & 48597 s& 17843 s& 32300 s& 98740 s&
		1.7 GB\\
		20 & 95457 s& 30467 s& 59384 s& 185309 s&
		2.8 GB\\
		21 & 180954s & 57621s & 126000s & 364576s
		& 4.7 GB \\
		\hline	
	\end{tabular}
	\caption{Timings and size of the calculated moments $m_l(n,k)$.}\label{Fig:Timing}
\end{figure}

We note that for $l\in\NN$ a closed form expression for $m_l(n,k)$ can be given in terms of the harmonic numbers $S_r(n)$ and $S_r(k)$ (resp.\ $\bar{S}_r(n)$ and $\bar{S}_r(k)$) with $1\leq r\leq l$. We note further that the sequences generated by the harmonic numbers $\{S_r(n)\mid r\geq1\}$ (resp. $\{\bar{S}_r(n)\mid r\geq1\}$) are algebraically independent among the rational sequences, i.e., the representation of $m_l(n,k)$ in terms of the harmonic numbers is optimal. Interestingly enough, the algebraic independence can be shown with the help of the summation paradigm of parameterized (creative) telescoping; see~\cite[Example~6.3]{Schneider:10c}. For more general classes of harmonic sums we refer also to~\cite{AS:18}.

\section{Asymptotic expansions}

Two days after Doron Zeilberger's first email, I got a slight update of his challenge:

\medskip

\begin{verbatim}
Finally, to get the donation in your honor you have to complete
the challenge. Find asymptotic expressions for $m_l(n)$ for
at least l=7 (it would be nice to also have l=8), and
prove (hopefully automatically) that
\end{verbatim}

\vspace*{-0.4cm}

$$\tt\lim_{n\to\infty} \frac{m_l(n)}{m_2(n)^{l/2}}$$ 

\vspace*{0.1cm}

\begin{verbatim}
equals 0 if l is odd and (2*l-1)(2*l-3)...1  if l is even.
This would be a partial elementary reproof of the
asymptotic normality proved for arbitrary r using "fancy" probability
of the asymptotic normality, in the Henze-Last paper.
\end{verbatim}

\medskip

Given the calculations of $m_l(n,k)$ we proceed as follows. For $l=2$ we enter the computed closed-form expression

\begin{mma}
\In m2=\frac{
	2 k
	-k^2
	-2 n
	-2 k n
	+2 k^2 n
	+2 n^2
	-k n^2
}{(n-1) n^2}\newline
\hspace*{1cm}-\frac{k (2+n) S[1,k]}{n}
+\frac{k (2+n) S[1,k]}{n}
+k^2 S[2,k]
-k^2 S[2,n];\\
\end{mma}

\medskip

\noindent and use the \texttt{HarmonicSums} command to calculate the first terms of the expansion in $n$:

\begin{mma}
\In SExpansion[m2, n, 3] /. { 
	LG[n] \to Log[n] + \gamma}\\
\Out \frac{7 (k-1) k}{6 n^3}
+\frac{k (-1+6 k)}{12 n^2}
-k^2S[2,\infty]
+\frac{\frac{k}{2}
	+k^2
	-2 k S[1,k]
}{n}
-k S[1,k]
+k^2 S[2,k]
+\big(
k
+\frac{2 k}{n}
\big) (\gamma +Log[n])
\\
\end{mma}
\medskip

\noindent where $\gamma\approx0.577216$ is Euler's constant, $S[2,\infty]=S_2(\infty)=\zeta(2)=\frac{\pi^2}{6}$, and $\zeta(z)=\sum_{i=1}^{\infty}\frac{1}{i^z}$ denotes the Riemann-Zeta function. Summarizing, we have calculated
\begin{equation}\label{Equ:M2Epansion}
\begin{split}
m_2(n,k)=&\frac{7 (k-1) k}{6 n^3}
+\frac{k (-25+18 k)}{12 n^2}
+\frac{k (-1+2 k)}{2 n}
-\big(
k
+\frac{2 k}{n}
\big) S_1({k})
\\
&+k^2 S_2({k})-k^2 \zeta(2)
+\big(
k
+\frac{2 k}{n}
\big) (\gamma +\log (n))+O(\tfrac1{n^4}).
\end{split}
\end{equation}
Analogously, we obtain, e.g., the expansion of $m_3(n,k)$:
\begin{align*}
m_3(n,k)&=-\frac{k \big(
	31-38 k+8 k^2\big)}{4 n^3}
-\frac{k \big(
	73-90 k+12 k^2\big)}{12 n^2}+\frac{\frac{1}{2} k (-1+6 k)
	-3 k (-1+2 k) \zeta(2)
}{n}\\
&+(\gamma +\log (n)) \big(
k
-\frac{k (-13+12 k)}{2 n^3}
-\frac{3 k (-3+2 k)}{n^2}
+\frac{6 k}{n}
+\frac{6 k S_1({k})}{n}
\big)\\
&
+\big(
-k
+\frac{k (-13+12 k)}{2 n^3}
+\frac{3 k (-3+2 k)}{n^2}
-\frac{6 k}{n}
\big) S_1({k})\\
&-\frac{3 k S_1({k})^2}{n}
+\big(
3 k^2
+\frac{3 k (-1+2 k)}{n}
\big) S_2({k})
-2 k^3 S_3({k})
-3 k^2 \zeta(2)
+2 k^3 \zeta(3)\\
&-\frac{3 k (\gamma +\log (n))^2}{n}
+\big(
-\frac{3 (k-2)^2 k}{n^3}
+\frac{3 (k-2) k}{n^2}
\big) \delta(k-2)+O(\tfrac1{n^4}).
\end{align*}
In total we calculated these expansions of $m_l(n,k)$ up to $l\leq 18$ with the following timings
\small
$$
\begin{array}{c||c|c|c|c|c|c|c|c|c|c|c|c|c|c}
l&5&6&7&8&9&10&11&12&13&14&15&16&17&18\\
\hline
\text{time}& 1.2s&3.6s&9.6s&25&65s&167s&405s&986s&2419s&1.6h&3.8h&9.8h&22.7h&53.4h.
\end{array}
$$
\normalsize
Here we parallelized the calculations on 10 Mathematica kernels. E.g., we needed $53.4$ hours to obtain the expansion for $m_{18}(n,k)$; however summing up all the timings of the used kernels we needed in total $16$ days of CPU time to tackle the case $l=18$. The expansions up to order 16 are online available and the corresponding link can be found on page~\pageref{Page:Links}.

To fully win Doron Zeilberger's challenge, we also dealt with the limit
$$c_l(k)=\lim_{n\to\infty} \frac{m_l(n,k)}{m_2(n,k)^{l/2}},\quad\quad l\geq2.$$
Looking at~\eqref{Equ:M2Epansion} the leading term in the expansion of $m_2(n,k)$ equals $\log(n)\,k$. In addition the leading term in $m_3(n,k)$ is also $\log(n)\,k$. Thus
$$c_3(k)=\lim_{n\to\infty} \frac{m_3(n,k)}{m_3(n,k)^{3/2}}=\lim_{n\to\infty}\frac{1}{\sqrt{k\log(n)}}=0.$$
Similarly, given the other computed expansions we get
$$c_{2l+1}=0,\quad 0\leq l\leq 8$$
and
$$\begin{array}{c||c|c|c|c|c|c|c|c|c}
l & 2& 4& 6& 8& 10& 12& 14& 16&18\\
\hline
c_l(k)& 1 & 3 & 15 & 105 & 945 & 10395 & 135135 & 2027025 & 34459425.\\
\end{array}$$
As a consequence we can confirm that the first values $c_{2l}(k)$ agree with the double factorial of odd numbers (sequence A001147 in OEIS), i.e.,
$$c_{2l}(k)=(2l-1)!!=\prod_{i=1}^{l}(2 i - 1)=\frac{(2 l - 1)!}{(l - 1)!2^{l-1}}$$
holds for all $1\leq l\leq 9$.

\section{Conclusion}

When I received Doron Zeilberger's email, the first calculations were straightforward -- except that I was not (and I am still not) expert in probability theory and thus made some annoying errors in the beginning. However, pushing the calculations further to the $21^{st}$ moment and computing its asymptotic expansions was a challenge. Here the computer algebra packages developed at RISC (\texttt{Sigma.m} and \texttt{EvaluateMultiSums.m} by myself and \texttt{HarmonicSums.m} by Jakob Ablinger) were of great help to calculate these huge expressions. I hope that this little note advertises how existing computer algebra tools in general and in particular those of the combinatorics group at RISC can push forward interesting (combinatorial) research topics. For instance, these packages have been heavily used in the last years to carry out large-scale QCD-calculations in the research field of elementary particle physics; see, e.g.,~\cite{CALadder:16} and references therein. Last but not least, I am very grateful to Doron Zeilberger who challenged me with these problems and initiated this fun project.

\medskip

\noindent\textbf{Acknowledgement.} I would like to thank the two referees for their valuable suggestions to improve the presentation.


\end{document}